\documentclass[a4paper]{amsart}
\begin{document}

\newcommand{\ea}{\mbox{{\bf a}}}
\newcommand{\eu}{\mbox{{\bf u}}}
\newcommand{\ueu}{\underline{\eu}}
\newcommand{\ueo}{\overline{u}}
\newcommand{\oeu}{\overline{\eu}}
\newcommand{\ew}{\mbox{{\bf w}}}
\newcommand{\ef}{\mbox{{\bf f}}}
\newcommand{\eF}{\mbox{{\bf F}}}
\newcommand{\eC}{\mbox{{\bf C}}}
\newcommand{\en}{\mbox{{\bf n}}}
\newcommand{\eT}{\mbox{{\bf T}}}
\newcommand{\eL}{\mbox{{\bf L}}}
\newcommand{\eR}{\mbox{{\bf R}}}
\newcommand{\eV}{\mbox{{\bf V}}}
\newcommand{\eU}{\mbox{{\bf U}}}
\newcommand{\ev}{\mbox{{\bf v}}}
\newcommand{\eve}{\mbox{{\bf e}}}
\newcommand{\uev}{\underline{\ev}}
\newcommand{\eY}{\mbox{{\bf Y}}}
\newcommand{\eK}{\mbox{{\bf K}}}
\newcommand{\eP}{\mbox{{\bf P}}}
\newcommand{\eS}{\mbox{{\bf S}}}
\newcommand{\eJ}{\mbox{{\bf J}}}
\newcommand{\eB}{\mbox{{\bf B}}}
\newcommand{\eH}{\mbox{{\bf H}}}
\newcommand{\leb}{\mathcal{ L}^{n}}
\newcommand{\eI}{\mathcal{ I}}
\newcommand{\eE}{\mathcal{ E}}
\newcommand{\hen}{\mathcal{H}^{n-1}}
\newcommand{\eBV}{\mbox{{\bf BV}}}
\newcommand{\eA}{\mbox{{\bf A}}}
\newcommand{\eSBV}{\mbox{{\bf SBV}}}
\newcommand{\eBD}{\mbox{{\bf BD}}}
\newcommand{\eSBD}{\mbox{{\bf SBD}}}
\newcommand{\ecs}{\mbox{{\bf X}}}
\newcommand{\eg}{\mbox{{\bf g}}}
\newcommand{\paromega}{\partial \Omega}
\newcommand{\gau}{\Gamma_{u}}
\newcommand{\gaf}{\Gamma_{f}}
\newcommand{\sig}{{\bf \sigma}}
\newcommand{\gac}{\Gamma_{\mbox{{\bf c}}}}
\newcommand{\deu}{\dot{\eu}}
\newcommand{\dueu}{\underline{\deu}}
\newcommand{\dev}{\dot{\ev}}
\newcommand{\duev}{\underline{\dev}}
\newcommand{\weak}{\rightharpoonup}
\newcommand{\weakdown}{\rightharpoondown}
\renewcommand{\contentsname}{ }

\newcommand{\down}[1]{{#1}^{\downarrow}}

\newtheorem{rema}{Remark}[section]
\newtheorem{thm}{Theorem}[section]
\newtheorem{lema}{Lemma}[section]
\newtheorem{prop}{Proposition}[section]
\newtheorem{defi}{Definition}[section]
\newtheorem{conje}{Conjecture}[]
\newtheorem{exempl}{Example}[section]
\newtheorem{opp}{Open Problem}[]
\renewcommand{\contentsname}{ }
\newenvironment{rk}{\begin{rema}  \em}{\end{rema}}
\newenvironment{exemplu}{\begin{exempl}  \em}{\end{exempl}}

\title{Lower semi-continuity of integrals with $G$-quasiconvex
potential}
\author{Marius Buliga}

\begin{abstract}
This paper introduces the proper notion of variational quasiconvexity associated 
to a group of diffeomorphisms. 
We prove a lower semicontinuity theorem connected to this notion. 
In the second part of the paper we apply this result to a 
class of  functions, introduced in \cite{buliga}. Such functions 
are $GL(n,R)^{+}$ quasiconvex, hence they induce lower semicontinuous integrals. 
\end{abstract}


\maketitle


\section{Introduction}

Lower semi-continuity of variational integrals
$$u \ \mapsto \ I(u) = \int_{\Omega} w(Du(x)) \mbox{ d}x$$
defined over Sobolev spaces is connected to the convexity of the 
potential $w$. In the scalar case, that is for functions $u$ with 
domain or range in $R$, the functional $I$ is weakly $W^{1,p}$ 
lower semi-continuous (weakly * $W^{1,\infty}$) if  and only if 
$w$ is convex, provided it is continuous and satisfies some 
growth conditions. The notion which replaces convexity in the vector 
case is quasi-convexity (introduced by Morrey \cite{9}). 

We shall concentrate on the case $u: \Omega \subset R^{n} \rightarrow R^{n}$
which is interesting for continuum media mechanics. Standard notation will be used, like:

\vspace{.3cm}

\begin{tabular}{ll}
$gl(n,R)$ & the linear space (Lie algebra) of $n\times n$ real matrices \\
$GL(n,R)$ & the group of invertible $n \times n$ real matrices \\
$GL(n,R)^{+}$ & the group of matrices with positive determinant \\ 
$sl(n,R)$ & the algebra of traceless $n\times n$ real matrices \\ 
$SL(n,R)$ & the group of real matrices with determinant one \\
$CO(n)$ & the group of conformal matrices \\ 
$id$ & the identity map \\
$1$ & the identity matrix \\
$\circ$ & function composition
\end{tabular}

\vspace{.3cm}

In this frame Morrey's  
quasiconvexity has the following definition. 

\begin{defi}
Let $\Omega \subset R^{n}$ be an open bounded set such that $\mid \partial \Omega 
\mid = 0$ and $w: gl(n,R) \rightarrow R$ be a measurable function. The map 
$w$ is quasiconvex if for any $H \in gl(n,R)$ and any Lipschitz $\eta: \Omega 
\rightarrow R^{n}$,  such that $\eta(x) = 0$ on $\partial \Omega$, we have
\begin{equation}
\int_{\Omega} w(H) \ \leq \ \int_{\Omega} w(H + D\eta(x)) 
\label{morreyqc}
\end{equation}
\label{dmorreyqc}
\end{defi}
Translation and rescaling arguments show that the choice of $\Omega$ is  irrelevant in the above definition.

Any quasiconvex function $w$ is rank one convex. There are several ways to define 
rank one convexity but this is due to the regularity assumptions upon $w$. The most 
natural, physically meaningful and historically justified, is to suppose that $w$ is $C^{2}$ and link rank one convexity with the ellipticity (cf. Hadamard 
\cite{11}) of the Euler-Lagrange 
equation associated to $w$. There are well-known ways to show that one can get rid 
of any regularity assumption upon $w$, replacing it by some growth conditions. Rank 
one convexity becomes then just what the denomination means, that is convexity along 
any rank one direction. 

\begin{prop}
Suppose that $w: gl(n,R) \rightarrow R$ is $C^{2}$ and quasiconvex. Then for any 
pair $a,b \in R^{n}$ the ellipticity inequality 
\begin{equation}
\frac{\partial^{2} w}{\partial H_{ij} \partial H_{kl}}(H) a_{i}b_{j}a_{k}b_{l} \ 
\geq \ 0
\label{ellip}
\end{equation}
holds true. 
\label{pellip}
\end{prop}

\begin{proof}
Take any $\eta \in C^{2}(\Omega, R^{n})$ such that $\eta(x) = 0$ on $\partial \Omega$ and $H \in gl(n,R)$. If $w$ is quasiconvex then the function 
$$t \mapsto \ f(t) \ = \ \int_{\Omega} w(H + t D\eta(x))$$
has a minimum in $t = 0$. Therefore $f'(0)=0$ and $f''(0)\geq 0$. Straightforward computation shows that $f'(0)=0$ anyway and $f''(0)\geq 0$ reads: 
$$\frac{\partial^{2} w}{\partial H_{ij} \partial H_{kl}}(H) \ \int_{\Omega} 
\eta_{i,j}(x)\eta_{k,l}(x) \ \geq \ 0$$
With the notation 
$$\Delta(\eta) \ = \ \int_{\Omega} D\eta(x) \otimes D\eta(x)$$
remark that $\Delta(x) \in V=gl(n,R) \otimes gl(n,R)$, because $V$ is a vectorspace 
and $D\eta(x) \otimes D\eta(x) \in V$ for any $x \in \Omega$. It follows that 
there is $P \in gl(n,R)$ such that: 
$$\Delta(\eta)_{ijkl} = P_{ij}P_{kl}$$
Integration by parts shows that $\Delta(\eta)$ has more symmetry, namely: 
$$\Delta(\eta)_{ijkl} = \Delta(\eta)_{ilkj}$$
which turns to be equivalent to $rank \ P \ \leq 1$. Therefore there are 
$a,b \in R^{n}$ such that $P = a\otimes b$. 

All it has been left to prove is that for any $a,b \in R^{n}$ there is  a $\lambda \not = 0$ and a  vector field $\eta \in C^{2}(\Omega, R^{n})$ such that $\eta(x) = 0$ on $\partial \Omega$ and $\Delta(\eta) = \lambda a \otimes b$. 
For this suppose that $\Omega$ is the unit ball in $R^{n}$, take $u: [0,\infty] \rightarrow R$ a $C^{\infty}$ map,  such that $u(1)=0$ and define: 
$$\eta(x) = u(\mid x \mid^{2}) \sin (b \cdot x) a$$
It is a matter of computation to see that  $\eta$ is well chosen to prove the thesis. 
\end{proof}

In elasticity  the elastic potential function $w$ is not defined on the Lie 
algebra $gl(n,R)$ but on the Lie group $GL(n,R)$ or a subgroup of it. It would be therefore 
interesting to find the connections between lower semicontinuity of the functional 
and the (well chosen notion of) quasiconvexity in this non-linear context. This is a 
problem which floats in the air for a long time. Let us recall two different definitions of quasiconvexity which are relevant. 
\begin{defi}
Let $w: GL(n,R)^{+} \rightarrow R$. Then: 
\begin{enumerate}
\item[(a)](Ball \cite{1}) $w$ is quasiconvex if for any $F \in GL(n,R)^{+}$ and 
any $\eta \in C^{\infty}_{c}(\Omega,R^{n})$ such that $F + D\eta(x) \in GL(n,R)^{+}$
for almost any $x \in \Omega$ we have 
$$\int_{\Omega} w(F + D\eta(x)) \ \geq \ \mid \Omega\mid w(F)$$
\item[(b)](Giaquinta, Modica \& Soucek \cite{gms}, page 174, definition 3) $w$ is 
Diff-quasiconvex if for any diffeomorphism $\phi: \Omega \rightarrow \phi(\Omega)$ 
such that $\phi(x)= Fx$ on $\partial \Omega$, for some $F \in GL(n,R)^{+}$ we have: 
$$\int_{\Omega} w(D\phi(x)) \ \geq \ \int_{\Omega} w(F)$$
\end{enumerate}
\label{dballgms}
\end{defi}
These two definitions are equivalent.

It turns out that very little is known about the lower semicontinuity properties 
of integrals given by Diff-quasiconvex potentials. It is straightforward that 
Diff-quasiconvexity is a necessary condition for weakly * $W^{1,\infty}$ (or uniform convergence of Lipschitz mappings) (see \cite{gms} proposition 2, same page). All that is known  reduces to the properties of polyconvex maps. A polyconvex map 
$w: GL(n,R)^{+} \rightarrow R$ is described by a convex function 
$g: D \subset R^{M} \rightarrow R$ (the domain of definition $D$ is convex as well) 
and $M$ rank one affine functions $\nu_{1},..., \nu_{M}: GL(n,R)^{+} \rightarrow R$
such that for any $F \in GL(n,R)^{+}$ 
$$w(F) = g(\nu_{1}(F),...,\nu_{M}(F))$$
The rank one affine functions are known(cf. Edelen \cite{5}, Ericksen \cite{6}, Ball, Curie, Olver \cite{21}): $\nu$ is rank one affine if and only if 
$\nu(F)$ can be expressed as a linear combination of subdeterminants of $F$ (uniformly with respect to $F$). Any rank one convex function is also called a null Lagrangian, because it generates a trivial Euler-Lagrange equation. 

Polyconvex function give lower semicontinuous 
functionals, as a consequence of Jensen's inequality and continuity of (integrals of) null lagrangians. This is a very interesting path to follow (cf. Ball \cite{1'})  and it leads to 
many applications. But it leaves unsolved the problem: are the integrals given 
by Diff-quasiconvex potentials lower semicontinuous? 

In the case of incompressible elasticity one has to work with the group of matrices 
with determinant one, i.e. $SL(n,R)$. The "linear" way of thinking has been compensated by wonders of analytical ingenuity. One purpose of this paper is to show 
how a slight modification of thinking, from linear to nonlinear, may give interesting 
results in the case $w: G \rightarrow R$ where $G$ is a Lie subgroup 
of $GL(n,R)$. Note that when $n$ is even a group which deserves attention is $Sp(n,R)$, the group of symplectic matrices. 

From now on linear transformations of $R^{n}$ and their matrices are identified. 
$G$ is a Lie subgroup of $GL(n,R)$.
\begin{defi}
For any $\Omega \subset
R^{n}$ open, bounded, with smooth boundary, we introduce the set 
$[G]^{\infty}(\Omega)$ of all bi-Lipschitz mappings $u$ from $\Omega$ 
to $R^{n}$ such that for almost any $x \in \Omega$ we have 
$Du(x) \in G$. 
\end{defi}
The set $Q \subset R^{n}$ is the unit cube $(0,1)^{n}$. 

The departure point of the paper is the following natural definition. 

\begin{defi}
The continuous function $w: G \rightarrow R$ is $G$-quasiconvex if 
for any $F \in G$ and $u \in [G]^{\infty}(Q)$ we have: 
\begin{equation}
\int_{Q} w(F) \mbox{ d}x \ \leq \ \int_{Q} w(FDu(x)) \mbox{ d}x
\label{edgqc}
\end{equation}
\label{dgqc}
\end{defi}

We describe now the structure of the paper. After the formulation of the lower 
semicontinuity theorem \ref{mainthm}, in section 3 is shown that quasiconvexity 
in the sense of definition \ref{dballgms} is the same as $GL(r,n)^{+}$ quasiconvexity. 
Theorem \ref{mainthm} is proved in section 4; in the next section is described 
the rank one convexity (or ellipticity) notion associated to $G$ quasiconvexity. 
The cases $GL(n,R)$ and $SL(n,R)$ are examined in detail. It turns out that classification of all 
universal conservation laws in incompressible elasticity is 
based on some unproved assumptions. In section 6 is described a class of 
$GL(n,R)^{+}$ quasiconvex functions introduced 
in Buliga \cite{buliga}. Theorem \ref{mainthm} is used to prove that any such function 
induces a lower semicontinuous integral.

\section{G-quasiconvexity and the lower semicontinuity result}
We denote by $[G]^{\infty}_{c}$  the class of all 
Lipschitz mapping from $R^{n}$ to $R^{n}$ such that 
$u - \ id$ has compact support and for almost any $x \in R^{n}$ we 
have $Du(x) \in G$. The main result of the paper is: 

\begin{thm}
Let $G$ be a Lie subgroup of $GL(n,R)$,  $\Omega$ an open,
bounded set with $\mid \partial \Omega \mid = 0$ and $w:G \rightarrow R$ 
locally Lipschitz. 
\begin{enumerate}
\item[a)] Suppose that for any sequence 
$u_{h} \in [G]^{\infty}_{c}$ weakly * $W^{1,\infty}$ convergent 
to $id$ we have: 
\begin{equation}
\int_{\Omega} w(F) \mbox{ d}x \ \leq \ \liminf_{h \rightarrow
\infty}  \int_{\Omega} w(F Du_{h}(x)) \mbox{ d}x
\label{agqc}
\end{equation}
Then for any bi-Lipschitz $u \in [G]^{\infty}_{c}$ and for any
sequence $u_{h}$ weakly * $W^{1,\infty}$ convergent 
to $u$ we have: 
\begin{equation}
\int_{\Omega} w(Du(x)) \mbox{ d}x \ \leq \ \liminf_{h \rightarrow
\infty}  \int_{\Omega} w(Du_{h}(x)) \mbox{ d}x
\label{lsc}
\end{equation}
Moreover, if \eqref{lsc} holds for any bi-Lipschitz 
$u \in [G]^{\infty}_{c}$ and for any
sequence $u_{h}$ weakly * $W^{1,\infty}$ convergent 
to $u$ then $w$ is $G$-quasiconvex. 
\item[b)] Suppose that $G$ contains the group $CO(R^{n})$ of
conformal matrices. Then \eqref{lsc} holds for any bi-Lipschitz 
$u \in [G]^{\infty}_{c}$ and for any
sequence $u_{h}$ weakly * $W^{1,\infty}$ convergent 
to $u$ if and only if $w$ is $G$-quasiconvex.
\end{enumerate}
\label{mainthm}
\end{thm}

The fact that weakly * lower semicontinuity implies $G$ quasiconvexity (end of 
point (a)) is 
easy to prove by rescaling arguments (cf. proposition 2,  Giaquinta, Modica and Soucek {\it op. cit.}). 

The method of proving the point (a) of the theorem is well known (see Meyers 
\cite{3}). Even if there is  nothing new there from the pure analytical viewpoint, I think that the proof deserves attention.

\section{G-quasiconvexity}

This section contains preliminary properties of $G$-quasiconvex 
continuous functions.

\begin{prop}
\begin{enumerate}
\item[a)] In the definition of $G$-quasiconvexity the cube 
$Q$ can be replaced by any open bounded set $\Omega$ such that 
$\mid \paromega \mid = 0$. 
\item[b)] The function $w$ is $G$-quasiconvex if and only if 
for any $F \in G$ and $u \in [G]^{\infty}_{c}(Q)$ we have: 
\begin{equation}
\int_{Q} w(F) \mbox{ d}x \ \leq \ \int_{Q} w(Du(x) F) \mbox{ d}x
\end{equation}
The converse is true. 
\item[c)] For any $U \in GL_{n}$ such that $ U G U^{-1} \subset G$ 
and for any $W: G \rightarrow R$  $G$-quasiconvex, the mapping  
$W_{U}: G \rightarrow R$, $W_{U} (F) =  W(U  F U^{-1}) $ 
is $G$-quasi-convex. 
\end{enumerate}
\label{prop1}
\end{prop}

\begin{rk}
The point b) shows that the non-commutativity of the multiplication 
operation does not affect the definition of $G$-quasiconvexity. 
 The point 
c) is a simple consequence of the fact that $G$ is a group. 
\end{rk}

\begin{proof}
The point a) has  a straightforward proof by translation and
rescaling arguments. 

For b) let us consider $F \in G$ and  an
arbitrary open bounded $\Omega \subset R^{n}$ with smooth boundary. 
The application which maps  $\phi \in [G]^{\infty}_{c}(\Omega)$ to 
 $F^{-1} \phi F \in 
[G]^{\infty}_{c}(F^{-1}(\Omega))$ is well defined and bijective. 
By a), if the function $w$ is $G$-quasiconvex then we have 
$$\int_{F^{-1}(\Omega)} w(F D(F^{-1} \phi F)(x)) \mbox{ d}x \geq 
\mid F^{-1}(\Omega) \mid w(F)$$
The change of variables $x = F^{-1}y$ resumes the proof of b). 

With $U$ like in the hypothesis of c), the application which maps 
$\phi \in [G]^{\infty}_{c}(\Omega)$ to $U \phi U^{-1} \in 
[G]^{\infty}_{c}(U^{-1}(\Omega))$ is well defined and bijective. The
proof resumes as for the point b). 
\end{proof}

The following proposition shows that quasi-convexity in the sense of definition 
\ref{dballgms} is a particular case of $G$-quasiconvexity.

\begin{prop}
Let us consider $F \in GL(n,R)^{+}$. Then  
$w$ is $GL(n,R)^{+}$-quasiconvex  in $F$  if
 and only if it is quasi-convex in $F$ in the sense of Ball.
\label{p4}
\end{prop}

\begin{proof}
Let  $E \subset R^{n}$ be an open bounded
 set and $\phi \in [GL(n,R)^{+}]^{\infty}_{c}(E)$. The vector field 
$\eta \ = \ F (\phi - id)$ verifies the condition that almost everywhere 
$F + D\eta(x)$ is invertible. Therefore, if 
$w$ is quasi-convex in  $F$, we derive from the inequality: 
$$\int_{E}w(F D\phi(y)) \mbox{ d}
 y \ \geq \ \mid E \mid 
W(F) \ \ . $$ 
We implicitly used the chain of equalities $$F + D 
\eta(y) \ =  \ F + F D\phi (y) \ - \  F \ =  
\ F D \phi (y) \ \ .$$ We have proved that quasi-convexity implies 
$GL(n,R)^{+}$-quasiconvexity. 

In order to prove the inverse implication let us  consider  
$\eta$ such that almost everywhere 
$F + D\eta(x)$ is invertible. We have therefore  
$\phi  \ = \ F^{-1} \psi \in [GL(n,R)^{+}](E)$ and 
$F D \phi \ = \ F + D \eta $. We use now the hypothesis that 
 $w$ is $GL(n,R)^{+}$-quasiconvex in $F$ and we find that $w$ is also 
quasi-convex. 
\end{proof}

\section{Proof of Theorem \ref{mainthm}}

The proof is divided into three steps. In the first step we shall 
prove the following: 

{\bf (Step 1.)}{\it Let $w: GL(n,R) \rightarrow R$  be locally Lipschitz. 
Suppose that for any Lipschitz bounded sequence $u_{h} \in [GL(n,R)]_{c}^{\infty}$ uniformly  convergent 
to $id$ on $\overline{\Omega}$ and for any  $F \in GL(n,R)$ we have: 
\begin{equation}
\int_{\Omega} w(F) \mbox{ d}x \ \leq \ \liminf_{h \rightarrow
\infty}  \int_{\Omega} w(F Du_{h}(x)) \mbox{ d}x
\label{1agqc}
\end{equation}
Then for any bi-Lipschitz $u: R^{n} \rightarrow R^{n}$ and for any sequence 
$u_{h}\in [GL(n,R)]_{c}^{\infty}$ uniformly convergent 
to $id$ on $\overline{\Omega}$ we have: 
\begin{equation}
\int_{\Omega} w(Du(x)) \mbox{ d}x \ \leq \ \liminf_{h \rightarrow
\infty}  \int_{\Omega} w(D(u_{h} \circ u)(x)) \mbox{ d}x
\label{llsc}
\end{equation}
}
\begin{rk}
This is just the point a) of the main theorem for the whole group of 
linear invertible transformations.
\end{rk}

\begin{proof}

For $\varepsilon > 0$ sufficiently small consider the set: 
$$U^{\varepsilon} \ = \ \left\{ B = \overline{B}(x,r) \subset \Omega \ : 
\ \exists \ A \in GL(n,R) \ , \ \int_{B} \mid Du(x) - A \mid < \varepsilon 
\mid B \mid \right\}$$
From the Vitali covering theorem and from the fact that $u$ is bi-Lipschitz we 
deduce that there is a sequence $B{j} = \overline{B}(x_{j}, r_{j}) \in 
U^{\varepsilon}$ such that: 
\begin{enumerate}
\item[-] $\mid \Omega \setminus \bigcup_{j}B_{j}\mid = 0$
\item[-] for any $j$ $u$ is approximatively differentiable in $x_{j}$ and 
$Du(x_{j}) \in GL(n,R)$
\item[-] we have
$$\int_{B_{j}} \mid Du(x) - Du(x_{j})\mid \ < \ \varepsilon \mid B_{j} \mid$$
\end{enumerate}
Choose $N$ such that 
$$\mid \Omega \setminus \bigcup_{j=1}^{N}B_{j} \mid \ < \ \varepsilon$$
We have therefore: 
$$\int_{\Omega}w(D(u_{h}\circ u)(x)) \ \geq \ \sum_{j=1}^{N} \int_{B_{j}} 
w(D(u_{h}\circ u)(x)) 
\ - \ C \varepsilon$$
$$\sum_{j=1}^{N} \int_{B_{j}} w(D(u_{h}\circ u)(x)) = J_{1} + J_{2} + J_{3}$$
where the quantities $J_{i}$ are given below, with their estimates. 
$$J_{1} = \sum_{j=1}^{N} \int_{B_{j}} \left[ w(Du_{h}(u(x))Du(x)) - 
w(Du_{h}(u(x))Du(x_{j}))\right]$$
$$\mid J_{1} \mid \leq \sum_{j=1}^{N} \int_{B_{j}} \mid w(Du_{h}(u(x))Du(x)) - 
w(Du_{h}(u(x))Du(x_{j}))\mid \ < \ C\varepsilon$$
$$J_{2} = \sum_{j=1}^{N} \int_{B_{j}} \left[ w(Du_{h}(u(x))Du(x_{j})) - 
w(Du_{h}(\overline{u}_{j}(x))Du(x_{j}))\right]$$
where $\bar{u}_{j}(x) = u(x_{j}) + Du(x_{j})(x-x_{j})$. We have the estimate: 
$$\mid J_{2} \mid \leq C \varepsilon$$
Indeed, by changes of variables we can write: 
$$I_{j}' \ = \ \int_{B_{j}}  w(Du_{h}(u(x))Du(x_{j})) \ = \ \int_{u(B_{j})}w(D u_{h}(y) Du(x_{j}) \ \mid \det Du^{-1}(y)\mid $$
$$I_{j}" \ = \ \int_{B_{j}} w(Du_{h}(\bar{u}_{j}(x))Du(x_{j})) \ = \ 
\int_{\bar{u}_{j}(B_{j})}w(D u_{h}(y) Du(x_{j}) \ \mid \det (Du(x_{j}))^{-1}\mid$$
The difference $\mid I_{j}' - I_{j}"\mid$ is majorised like this
$$ \mid I_{j}' - I_{j}"\mid \ \leq \ \int_{u(B_{j}) \cap \bar{u}_{j}(B_{j})} 
C  \mid \mid \det Du^{-1}(y)\mid - \mid\det (Du(x_{j}))^{-1}\mid \mid  \ + \ C \mid u(B_{j}) \Delta \bar{u}_{j}(B_{j}) \mid$$
The function $\mid \det \cdot \mid$ is rank one convex and satisfies the growth condition 
$\mid \det F \mid \ \leq c(1+\mid F \mid^{n})$ for any $F \in GL(n,R)$. Therefore this function satisfies also the inequality: 
$$\mid \mid \det F \mid - \mid \det P \mid \mid \ \leq \ C \mid F - P \mid 
\left( 1 + \mid F \mid^{n-1} + \mid P \mid^{n-1}\right)$$
Use now this inequality, the properties of the chosen Vitali covering and the 
uniform bound on Lipschitz norm of $u$, $u_{h}$, to get the claimed estimate.

$$J_{3} = \sum_{j=1}^{N} \int_{B_{j}} w(Du_{h}(\overline{u}_{j}(x))Du(x_{j}))$$
By the change of variable $y = \overline{u}_{j}(x)$ and the hypothesis we have 
$$\liminf_{h\rightarrow \infty}J_{3} \geq \ \liminf_{h\rightarrow \infty}
\sum_{j=1}^{N} \int_{B_{j}}w(Du(x_{j}))$$
Put all the estimates together and pass to the limit with $N \rightarrow 
\infty$ and then $\varepsilon \rightarrow 0$. 
\end{proof}

{\bf (Step 2.)}{ If we replace in }{\bf Step 1.} {\it the group 
$GL(n,R)$ by a Lie subgroup $G$ the conclusion is still true.}

\begin{proof}
Indeed, remark that in the proof of the previous step it is used only the fact 
that $GL(n,R)$ is a group of invertible maps. 
\end{proof}

{\bf Step 3.} {\it The point b) of the  Theorem \ref{mainthm} is
true.} 

\begin{rk}
In the classical setting of quasiconvexity, this step is proven by 
an argument involving Lipschitz extensions with controlled 
Lipschitz norm. In our case the corresponding 
Lipschitz extension assertion would be: 
{\it let $u \in [G]^{\infty}_{c}$ with Lipschitz norm 
$\| u - \ id \| = 
\varepsilon$. For $\delta > 0$ sufficiently big there exists 
$v \in [G](B(0, 1+ \delta))$ such that $v = u$ on $B(0,1)$ and 
$\| v - \ id\|$ controlled from above by $\varepsilon$.}  
This is not known to be true, even for $G = GL(n,R)$. That is why 
we shall use a different approach. 
\end{rk}

\begin{proof}
Because $G$ is a group, it is sufficient to make the proof for 
$F = 1$. 

Let $u_{h}\in [G]^{\infty}_{c}$ be a sequence weakly * convergent to 
$id$ on $\Omega$ and $D \subset \subset \Omega$. 
 For $\varepsilon > 0$  sufficiently small and $C>1$ we have
$$D_{C\varepsilon} \ = \ \bigcup_{x \in D} B(x,C\varepsilon) \ \subset \ 
\Omega$$
It is not restrictive to suppose that 
$$\lim_{h \rightarrow \infty} \int_{\Omega} w(Du_{h}) \mbox{ d}x $$ 
exists and it is finite. For any $\varepsilon > 0$ there is 
$N_{\varepsilon}$ such that for any $h > N_{\varepsilon}$ 
$u_{h}(D) \subset D_{\varepsilon}$.

Take a minimal Lipschitz extension 
$$\overline{u}_{h} : D_{C\varepsilon}\setminus C \rightarrow R^{n} \ , \ 
\overline{u}_{h}(x) = \left\{ \begin{array}{ll}
u_{h}(x) & , \ x \in \partial D \\ 
x & , \ x \in \partial D_{C\varepsilon}
\end{array}\right.$$
The Lipschitz norm of this extension, denoted by $k_{h}$,  is smaller than some constant independent on $h$. 

 Now, for any $h$ define: 
$$\psi_{h} \ = \ \frac{1}{2k_{h}} \ 
\overline{u}_{h_{|_{D_{C\varepsilon}} \setminus D}}$$
According to Dacorogna-Marcellini Theorem 7.28, Chapter 7.4. \cite{dacomarc}, there is a solution 
$\sigma_{h}$ of the problem 
$$\left\{ \begin{array}{ll}
D \sigma_{h} \in O(n) & \mbox{ a. e. in } D_{C\varepsilon} 
\setminus D \\
\sigma_{h} \ = \ \psi_{h} & \mbox{ on } \partial ( D_{\varepsilon} 
\setminus D 
\end{array} \right. $$
Let 
$$v_{h}(x) \ = \ \left\{ \begin{array}{ll}
u_{h}(x) & x \in D \\ 
k_{h} \sigma_{h}(x) & x \in \Omega \setminus D
\end{array} \right. $$
Note that $Dv_{h} \in CO(n)$. 

The following estimate is then true: 
$$\mid \int_{D} w(Du_{h}) \mbox{ d}x \ - \ \int_{\Omega} w(Dv_{h}) 
\mbox{ d}x \mid \ = \ \mid \int_{ D_{C\varepsilon} \setminus D} 
w(Dv_{h}) \mbox{ d}x \mid \ \leq \ $$
$$\leq \ \int_{ D_{C\varepsilon} \setminus D} 
\mid w(Dv_{h}) \mid \mbox{ d}x \ \leq \ C \mid D_{\varepsilon} 
\setminus D \mid$$ 
$w$ is $G$-quasiconvex, therefore: 
$$\int_{D_{\varepsilon}} w(Dv_{h}) \mbox{ d}x \ \geq \ \mid
D_{\varepsilon} \mid w(1)$$
We put all together and we get the inequality: 
$$\lim_{h \rightarrow \infty}\int_{D} w(Du_{h}) \mbox{ d}x \ \geq \ 
\mid D_{\varepsilon} \mid w(1) - C \mid D_{\varepsilon} 
\setminus D \mid$$ 
The proof finishes after we pass $\varepsilon$ to 0. 
\end{proof}

\section{Rank one convexity}

The rank-one convexity notion associated to $G$ quasi-convexity is described in the 
next proposition, for $w \in C^{2}(G,R)$. Before this, let us introduce a differential operator naturally connected to the group structure of 
$G$. Denote by $\mathcal{G}$ the Lie algebra of $G$. For any pair 
$(F,H) \in G \times \mathcal{G}$, the derivative of $w:G \rightarrow R$ in 
$F$ with respect to $H$ is 
$$Dw(F)H \ = \ \frac{d}{dt}_{|_{t=0}} w(F \exp(tH))$$
We shall also use the notation (for $F \in G$ and $H,P \in \mathcal{G}$): 
$$D^{2}w(F)(H,P) \ = \ D(Dw(\cdot)H)(F)P$$

\begin{prop}
A necessary condition for $w \in C^{2}(G,R)$ to be $G$ quasi-convex is 
$$\int_{\Omega} D^{2}w(F)(D\eta(x), D\eta(x)) \ = 0$$
for any $F \in G$ and $\eta \in C^{2}(\Omega,R^{n})$, $D\eta(x) \in  \mathcal{G}$ 
a.e. in $\Omega$, $supp \ \eta \in \Omega$. 
\label{proc}
\end{prop}

\begin{proof}
Given such an $\eta$, consider the solution of the o.d.e. problem: 
$$\dot{\phi}_{t} \ = \eta \circ \phi_{t} \ , \ \phi_{0} = \ id_{|_{\Omega}}$$
This is an one-parameter group in the diffeomorphism class 
$[G]^{\infty}(\Omega)$. Define then: 
$$f(t)  \ = \ \int_{\Omega} w(FD\phi_{t}(x))$$
The $G$ quasiconvexity of $w$ implies that $f$ has a  minimum in $t=0$. 
That means $f'(0) = 0$ and $f''(0) \geq 0$. The first condition is trivially 
satisfied and the second is, by straightforward computation, just the conclusion of 
the proposition. 
\end{proof}

We shall call $G$ rank one convex a function which satisfies the conclusion of 
the proposition \ref{proc}. 

Consider the vector space 
$$V(\mathcal{G}) = \left\{ (H,H) \in \mathcal{G} \times \mathcal{G} \ : \ H \in 
\mathcal{G} \right\}$$
and the set
$$RO(\mathcal{G}) = \left\{ (a,b) \in R^{n} \times R^{n} \ : \  
a \otimes b  \in \mathcal{G} \right\}$$

\begin{prop}
Suppose that $w:G \rightarrow R$ is a $C^{2}$ function. If for any $a,b \in RO(\mathcal{G})$ 
\begin{equation}
D^{2}w(F)(a\otimes b, a\otimes b) \geq 0
\label{efunny}
\end{equation}
then $w$ is $G$ rank one convex. 
\label{funnyprop}
\end{prop}

\begin{proof}
We shall use the notations from the proof of the preceding proposition. 
We see that 
$$\int_{\Omega} (D\eta(x), D\eta(x)) \ \in V(\mathcal{G})$$
Therefore there is an $X \in \mathcal{G}$ such that 
$$(X,X) = \int_{\Omega} (D\eta(x), D\eta(x))$$
Using integration by parts we find that for any indices 
$i,j,k,l \in 1, ... , n$ we have: 
$$X_{ij}X_{kl} = X_{il}X_{kj}$$
which implies that $X$ has rank one. Hence there are $a,b \in R^{n}$ such that 
$X= a\otimes b$. Use the definition of $G$ rank one convexity to prove that 
\eqref{efunny} implies the $G$ rank one convexity. 
\end{proof}

In the case $G=GL(n,R)$ we find that $GL(n,R)$ rank one convexity is equivalent to  classical 
rank one convexity.  To see 
this, take arbitrary $F \in GL(R^{n})$, $a,b \in R^{n}$, $s>0$ and 
$u \in C^{\infty}_{c}(\Omega,R)$. Define 
$$\eta^{s}(x) = u(x) \sin \left[s(b \cdot x)\right] \ a$$
Because $GL(n,R)$ is an open set in the vectorspace of $n\times n$ real matrices, 
the $GL(n,R)$ rank one condition reads: 
$$s^{2} \ \frac{d^{2} \ w}{d F_{ij} d F_{kl}}(F) (Fa)_{i} b_{j} (Fa)_{k} b_{l} 
 \int_{\Omega} u^{2}  \ + B \ \geq 0$$
with $B$ independent on $s$. We deduce that 
$$ \frac{d^{2} \ w}{d F_{ij} d F_{kl}}(F) (Fa)_{i} b_{j} (Fa)_{k} b_{l} \geq 0$$
for any choice of $F$, $a,b$. This is the same as: 
$$ \frac{d^{2} \ w}{d F_{ij} d F_{kl}}(F) a_{i} b_{j} (a_{k} b_{l} \geq 0$$
for any $F$, $a,b$. 

For the group $SL(n,R)$ of matrices with determinant one we obtain a similar 
condition by imposing the constraint $div \ \eta^{s} =0$. This can be done 
if $a \cdot b = 0$ and $Du(x) \cdot a = 0$. For simplicity suppose that 
$w$ is defined in a neighbourhood of $SL(n,R)$. Then $w$ is $SL(n,R)$ 
rank one convex implies 
\begin{equation}
 \frac{d^{2} \ w}{d F_{ij} d F_{kl}}(F) (Fa)_{i} b_{j} (Fa)_{k} b_{l} \geq 0
\label{slrankone}
\end{equation}
for any $F \in SL(n,R)$, $a,b \in R^{n}$, $a \cdot b = 0$.

\subsection{Rank one affine functions} 
A map $w$ is $G$ rank one affine if $w$ and $-w$ are $G$ rank one convex. 
For the case $G=GL(n)$ we see that the rank one affines are known. This 
is very useful in several instances. The reason is that the Euler-Lagrange 
equation associated to the potential $w$ does not change if one adds a rank 
affine function to $w$. At the action functional level 
$$I_{w}(\phi) \ = \ \int_{\Omega} w(D\phi(x))$$
the addition of a $GL(n,R)$ rank one function means the addition of a closed 
form which cancels with the integral. This coincidence led to the development 
of formal calculus of variations in the frame of the jet bundle formalism, which 
permits to classify all universal conservation 
laws in elasticity. For this classification  see Olver \cite{olver}. 

The case $G=SL(n,R)$ is equally important, because it is about incompressible 
elasticity. Or, in this case nothing is known, because it is not proven that 
the $SL(n,R)$ rank one affine functions correspond to closed forms. For this 
reason  Olver's classification \cite{olver} of universal conservation laws  
is not proven to be complete.  

We arrived to the following 

{\bf Open problem:} Describe all $G$ rank one affine functions. 

In particular situations the problem  has been solved. For example if 
$G= GL(n,R)$ then any rank one affine function is a classical null lagrangian.  
In the case $SL(2,R)$ we have the following theorem: 

\begin{thm}
Any $SL(2,R)$ rank one affine function is affine. 
\end{thm}

\begin{proof}
We have to prove that if $w: SL(2,R) \rightarrow R$ is rank one affine 
then $w(F) = a_{ij}F_{ij} + b$. It is sufficient to prove the thesis for 
any $F$ in an open dense set in $SL(2,R)$. We shall use the following 
maps: 
$$(X,Y,Z) \in R^{*} \times R \times R \ \mapsto \ 
F \ = \ \left( \begin{array}{cc}
X & Y \\
Z & \frac{1+YZ}{X}
\end{array} \right)$$
$$(X',Y',Z') \in R^{*} \times R \times R \ \mapsto \ 
F \ = \ \left( \begin{array}{cc}
\frac{1+Y'Z'}{X'} & Y' \\
Z' & X'
\end{array} \right)$$
Take arbitrary $a = (a_{1},a_{2})$ and perpendicular $b= (-a_{2},a_{2})$. 
If $w$ is $SL_{2}$ rank one affine then the mapping 
$$t \mapsto f(t;a\otimes b , F) =   w(F(1+ t a \otimes b))$$
is linear for any $F \in SL(2,R)$. We have used here the relation \label{slrankone} and the equality 
$\exp a \otimes b = 1+ a \otimes b$, for any orthogonal $a,b$. Rank one convexity 
of $w$ means that the second derivative of 
$f(t;a\otimes b , F)$ with respect to $t$ vanishes for any choice of $F$ and $a$. 

We express $F$ in terms of the coordinates $F=F(X,Y,Z)$ and $F=F(X',Y',Z')$. After some 
elementary computation we obtain the following minimal system of equations 
for the function $w(X,Y,Z) = w(F(X,Y,Z))$: 
\begin{equation}
\left\{ 
\begin{array}{lll}
w_{XX} X^{2} & = & 2 w_{YZ} (1+YZ) \\ 
w_{ZZ} X & = & - w_{YZ} Y \\ 
w_{XY} X & = & - w_{YZ} Z \\
w_{YY} & = & 0 \\ 
w_{ZZ} & = & 0 
\end{array} \right.
\label{sistem}
\end{equation}
From equations (\ref{sistem}.4) and (\ref{sistem}.5) we find that $w$ has the form 
$$w(X,Y,Z) = A(X)YZ + B(X)Y + C(X) Z + D(X)$$
From (\ref{sistem}.2) we obtain the equation 
$$XC'(X) + XY A'(X) = - A(X)Y$$
From here we derive that $C(X) = c$ and $A(X) = k/X$. We update the form of $w$, 
 use (\ref{sistem}.3) to get $B(X) = b$ and (\ref{sistem}.1) to get 
$D(X) = (k/X) + eX + f$. We collect all the information and we obtain that 
$w$ has the expression: 
$$w(X,Y,Z) = k \frac{1+YZ}{X} + bY + cZ + eX + f$$
which proves the theorem. 
\end{proof}

Therefore, in the case $G=SL(2,R)$ we have proved that there are no rank one 
affine functions other than the classical ones. The proof is not adapted 
to generalizations. The case $G= SL(3,R)$ is open.

Other groups are equally significant, like the group $Sp(n,R)$ of 
symplectomorphisms. I don't know of any attempt to solve this problem.

\subsection{Rank one convexity and quasiconvexity}
The $GL(n,R)$ rank one convexity is not equivalent to $GL(n,R)$ quasiconvexity 
in any dimension. 

\begin{prop}
The function $w: GL(n,R) \rightarrow R$ defined by 
$$w(F) \ = \ - \log \mid \det F \mid$$
is $GL(n,R)$ rank one convex but not $GL(n,R)$ quasiconvex. 
\label{pexe}
\end{prop}

\begin{proof}
The map is polyconvex hence it is rank one convex. It is not quasi-convex though. 
To see this fix $\varepsilon \in (0,1)$, $A \in GL(n,R)$ and $\Omega = B(0,1)$. There is a 
Lipschitz solution to the problem 
$$\left\{ \begin{array}{ll}
Dv(x) \in O(n) & \mbox{ a.e. in } \Omega \\
v(x) = \varepsilon x & x \in \partial \Omega 
\end{array} \right. $$
We have then, for $u(x) = v(x)/\varepsilon \ \in [GL(n,R)]^{\infty}(\Omega)$: 
$$\int_{\Omega} w(A Du(x)) \ = \ \int_{\Omega} - \log \mid \det A \mid \ + \ 
\int_{\Omega} n \log \varepsilon \ < \ \int_{\Omega} w(A) $$
\end{proof}

Next proposition justifies this result. 
\begin{prop}
For any $w: G \rightarrow R$ define $\i w: G \rightarrow R$ by: 
$$\i w (F) = \mid \det F \mid \ w(F^{-1})$$
Then $w$ is $G$ rank one convex if and only if $\i w$ is. Also, if $w$ is 
$G$ quasi-convex then for any $u \in [G]^{\infty}(\Omega)$ we have: 
$$\int_{\Omega} w(F Du(x)) \ \geq \int_{\Omega} w(F)$$
\label{puseful}
\end{prop}

\begin{proof}
Take $u$ like in the hypothesis. Then for any (continuous) $w$ we have 
$$\int_{\Omega} w(Du^{-1}(x)) \ = \ \int_{\Omega} \i w(Du(x))$$
by straightforward computation. Use now the proof of proposition \ref{proc}
to deduce the first part of the conclusion. For the second part use 
 the definition \ref{dgqc} and the proposition \ref{edgqc}. 
\end{proof}

Let us apply this proposition to $w(F) = - \log \mid \det F \mid$. Remark that 
when $\det F$ goes to zero the function goes to $+\infty$. Now, 
$\i w(F) \ = \ \mid \det F \mid \log \mid \det F \mid$ and this function can 
be continuously prolongated to matrices with determinant zero by setting 
$\i w(F) = 0$ if $\det F = 0$. It is easy to see that the prolongation 
of $\i w$ ceases to be rank one convex.

\section{Application: a class of quasiconvex functions}

The goal of this section is to give a class of quasi-convex isotropic functions which seem to be  
complementary to the polyconvex isotropic ones. We quote the following 
result of Thompson and Freede \cite{thomfre}, Ball \cite{1} (for a proof 
coherent with this paper see Le Dret \cite{ledret}). 

\begin{thm}
Let $g:[0,\infty)^{n} \rightarrow R$ be convex, symmetric and nondecreasing in each variable. 
Define the function $w$ by
$$w:gl(n,R) \rightarrow R \ , \ \ w(F) = g(\sigma(F)).$$
Then $w$ is convex. 
\label{thomfreede}
\end{thm}

We shall use the Theorem 6.2. Buliga \cite{buliga}. We need a
notation first. Let $x = (x_{1}, ... , x_{n}) \in R^{n}$ be a
vector. Then the vector $\down{x} = (\down{x}_{1}, ... ,
\down{x}_{n}) \in R^{n}$ is obtained by rearanging in decreasing
order the components of $x$. Remark that for any symmetric function 
$h : R^{n} \rightarrow R$ there exists and it is unique the function 
$p: R^{n} \rightarrow R$ defined by the relation: 
$$p(\sum_{i=1}^{k} \down{x}_{i}) \ = \ h(x_{k})$$

\begin{thm}
Let $g:(0,\infty)^{n} \rightarrow R$ be a continuous symmetric 
function and 
 $h:R^{n} \rightarrow R$, $h(x_{1}, ..., x_{n})  = g(\exp x_{1},
..., \exp x_{n})$. Suppose that 
\begin{enumerate}
\item[(a)] $h$ is convex, 
\item[(b)] The function $p$ associated to $h$  is 
nonincreasing in each argument. 
\end{enumerate}

  Let $\Omega \subset R^{n}$ be bounded, with piecewise 
smooth boundary and $\phi: \overline{\Omega} \rightarrow R$ be any Lipschitz function  such that $D\phi(x) \in GL(n,R)^{+}$ a.e. and $\phi(x) = x$ on $\partial \Omega$. 
Define the function
$$w: GL(n,R)^{+} \rightarrow R \ , \ \ w(F) = g(\sigma(F))$$
Then for any $F \in GL(n,R)^{+}$ we have: 
\begin{equation}
\int_{\Omega} w(F D\phi(x)) \ \geq \ \mid \Omega \mid w(F)
\label{firstqc}
\end{equation}
\label{texe}
\end{thm}

A consequence of  theorem \ref{texe} and Theorem \ref{mainthm} (a) is: 

\begin{prop}
In the hypothesis of Theorem \ref{texe}, let $\phi_{h}: \Omega \rightarrow R^{n}$ 
be a sequence of Lipschitz bounded functions such that 
\begin{enumerate}
\item[(a)] for any $h$ $D\phi_{h}(x) \in GL(n,R)^{+}$ a.e. in $\Omega$. 
\item[(b)] the sequence $\phi_{h}$ converges uniformly to $u: \Omega \rightarrow 
\Omega$, bi-Lipschitz function. 
\end{enumerate}
Then 
\begin{equation}
\liminf_{h \rightarrow \infty} \int_{\Omega} w(D\phi_{h}(x)) \ \geq \ \int_{\Omega} 
w(Du(x))
\label{firstlsc}
\end{equation}
\label{corlsc}
\end{prop}

\begin{proof}
It is clear that theorem \ref{texe} implies the hypothesis of point (a), theorem \ref{mainthm}. Indeed, the conclusion of theorem \ref{texe} can be written like this: 
for any  $u \in [GL(n,R)^{+}](\Omega)$ such that 
$$\bar{Du}(\Omega) = \frac{1}{\mid \Omega \mid} \int_{\Omega} Du(x) \mbox{ d}x \ \in GL(n,R)^{+}$$
we have the inequality 
$$\int_{\Omega} w(Du(x)) \mbox{ d}x \geq \ \int_{\Omega} w(\bar{Du}(\Omega)) \mbox{ d}x$$
Take a sequence of mapping $(u_{h}) \subset [GL(n,R)^{+}](\Omega)$ uniformly convergent 
to $F\in GL(n,R)^{+}$. The previous inequality and the continuity of $w$ imply: 
$$\int_{\Omega} w(F) \mbox{ d}x \ \leq \ \int_{\Omega}w(Du_{h}(x)) \mbox{ d}x$$
Apply now theorem \ref{mainthm} (a) and obtain the thesis.
\end{proof}

The class of functions $w$ described in theorem \ref{texe} and the class of 
polyconvex functions seem to be  different. However, by picking $h$ linear, we obtain 
a polyconvex function, like $$w(F) = - \log \mid \det F \mid$$
We have seen in proposition \ref{pexe} that this function is not $GL(n,R)$ quasiconvex but proposition \ref{corlsc} tells that $w$ is $GL(n,R)^{+}$ quasiconvex. 

We close with an example of another function which we can prove that it 
is $GL(n,R)^{+}$ quasiconvex. We use the notation 
$F = R_{F}U_{F}$ for the polar decomposition of $F \in GL(n,R)^{+}$, with 
$U_{F}$ symmetric and positive definite. The example is the function: 
$$w:GL(n,R)^{+} \rightarrow R \ , \ \ w(F) = \det F \ \log \left( trace \ U_{F} 
\right)$$
With the notation introduced in proposition \ref{puseful}, 
let's look to the  the function $\hat{w} = \i w$. It has the expression: 
$$\hat{w} : GL(n,R)^{+}  \rightarrow R \ , \ \ \hat{w}(F) \ = \ \log \left( \ 
trace \ U_{F}^{-1} \right)$$
It is a matter of straightforward computation to check that $\hat{w}$ verifies the 
hypothesis of theorem \ref{texe}. It is therefore $GL(n,R)^{+}$ quasiconvex. 
By proposition \ref{puseful} $w$ is $GL(n,R)^{+}$ quasiconvex, too,   hence 
lower semicontinuous in the sense of theorem \ref{mainthm} (a).

\vspace{2.cm}

\noindent
Institute of Mathematics of the Romanian Academy \\ 
 PO BOX 1-764, RO 70700, Bucharest, Romania, 
 e-mail: Marius.Buliga@imar.ro \\ and \\ 
Ecole Polytechnique F\'ed\'erale de Lausanne, D\'epartement de Math\'ematiques \\ 
1015 Lausanne, Switzerland, e-mail: Marius.Buliga@epfl.ch

\end{document}